\documentclass{article}
\usepackage[utf8]{inputenc}
\usepackage{amsmath}
\usepackage{amssymb}
\usepackage{amsfonts}
\usepackage{hyperref}

\title{On Power Sums and the Kummer Congruences}
\author{Samuel Goodman}
\date{May 2023}

\begin{document}
\maketitle
\section*{Introduction}
The simplest form of the Kummer Congruences were first proven in $1851$ by Ernst Eduard Kummer in [2]. They state that given a prime $p$ and even positive integers $r,s$ with $r$ not a multiple of $p-1$ such that $r \equiv s \mod p-1$, then $$\displaystyle \frac{B_r}{r} \equiv \frac{B_s}{s} \mod p$$ where $B_n$ denotes the $n$th Bernoulli number. A stronger version of this result was first obtained by Voronoi in [3] using power sum identities and then later by Kubota-Leopoldt in [1] using the theory of $p$-adic $L$-functions, proving that given a prime $p$, a nonnnegative integer $a$, and even positive integers $r,s$ with $r$ not a multiple of $p-1$ such that $r \equiv s \mod p^a(p-1)$, then $$\displaystyle (1-p^{r-1})\frac{B_r}{r} \equiv (1-p^{s-1})\frac{B_s}{s} \mod p^{a+1}$$ In particular, Kubota-Leopoldt cast this classical result in a more natural light, as it in fact becomes the input for the proof of the continuity of the $p$-adic zeta function. \newline \newline In this paper, we develop a structure theory that again recasts this result in a new way, shedding light on the element-wise behavior of a special class of polynomials $\mod$ powers of $p$ through a series of new results about the ``pseudo-multiplicative" structure these polynomials have by analyzing the stabilizers of the natural multiplication action on the resulting sets of values. As a result of this approach, we also deduce a new proof of the Kummer Congruences.
\newline \newline We start by introducing some preliminary results, well-known in the literature:
\subsection*{Lemma 1}
Let $p>3$ be a prime, $a$ a positive integer, and $r$ a positive even integer. We have that $$\displaystyle \sum _{n=1}^{p^a}n^r \equiv p^aB_r \mod p^{2a+v_p(r)+1}$$
\subsection*{Lemma 2}
Let $p$ be an odd prime and $a$ and $n$ be positive integers and $k$ a positive integer with $k \geq r+a$, $p^r|k$, and $p-1$ not dividing $k$. Then we have that $$\displaystyle \sum_{n=1}^{p^a} n^k \equiv 0 \mod p^{r+a}$$  
\subsection*{Lemma 3 (Adams Theorem)} For any even integer $s$ with $s$ not a multiple of $p-1$, $v_p(\frac{B_s}{s}) \geq 0$ (in particular, $\frac{B_s}{s}$ is a $p$-integer). \newline \newline We start by showing that certain special multisets are stabilized by nontrivial residue classes upon the natural multiplication action. This will be the first hint of local information (namely the individual values of the polynomial $\mod$ powers of $p$) have rich internal structure that explains global phenomena, one of which will be a precursor to the Kummer Congruences.
\section*{Theorem 1}
Let $p$ be an odd prime, $a,t$ nonnegative integers, and $k$ a positive integer with $p^{2a+1}|k$. Let $d=\frac{p-1}{\gcd(k,p-1)}$, $v=\min(v_p(k)-2a-1,t)$, and $g$ a $d$th root of unity $\mod p^{3a+t+v+2}$. In addition, let $\mathfrak{S}$ be the multiset consisting of the invertible residues $\mod p^{3a+t+v+2}$ attained by $n^{(k+p^a(p-1))p^t}+n^{(k-p^a(p-1))p^t}$ as $n$ varies from $1$ to $p^{a+1}$ with multiplicity. Then $g\mathfrak{S}=\mathfrak{S}$. \newline \newline
We show that $g\mathfrak{S} \subset \mathfrak{S}$, from which it follows that $g\mathfrak{S}=\mathfrak{S}$ as they are multisets of the same size. Each element $s$ of $\mathfrak{S}$ corresponds to some positive integer $n$ with $1 \leq n<p^{a+1}$ with $n^{(k+p^a(p-1))p^t}+n^{(k-p^a(p-1))p^t} \equiv s \mod p^{3a+t+v+2}$ by construction. For this $n$, we find an associated $n'$ such that $n'^{(k+p^a(p-1))p^t}+n'^{(k-p^a(p-1))p^t} \equiv gs \mod p^{3a+t+v+2}$. \newline \newline Let $k=p^rk'$ with $p \nmid k'$. First we show that there is a solution to $x^{k'} \equiv g \mod p^{a+1}$. Note that $\gcd(k',p-1)=\gcd(k,p-1)$, and so $x^{k'}$ and $x^{\gcd(k,p-1)}$ generate the same residue classes $\mod p^{a+1}$, so it suffices to show the result for $x^{\gcd(k,p-1)}$. Let $g'$ be a primitive root $\mod p^{a+1}$ and write $g \equiv g'^b \mod p^{a+1}$ for some $b$. Then since $g^d \equiv 1 \mod p^{a+1}$, we conclude that $p-1|db$, which implies that $\frac{p-1}{d}|b$, which gives that $\gcd(k,p-1)|b$ by definition and hence that $g \equiv x^{\gcd(p-1,k)} \mod p^{a+1}$ for some $x$. \newline \newline Now pick a fixed solution $x'$ to $x^{k'}\equiv g \mod p^{a+1}$. For each $n$, we let $n'$ be the unique positive integer $\leq p^{a+1}$ such that $n' \equiv nx' \mod p^{a+1}$, and we show that $n'^{(k+p^a(p-1))p^t}+n'^{(k-p^a(p-1))p^t} \equiv gs \mod p^{3a+t+v+2}$ under this assumption. Note that this is a bijection since $x'$ is coprime to $p$. \newline \newline
First note that $$(nx')^{(k+p^a(p-1))p^t}+(nx')^{(k-p^a(p-1))p^t}=(n'+p^{a+1}c)^{(k+p^a(p-1))p^t}+(n'+p^{a+1}c)^{(k-p^a(p-1))p^t}$$ but then the binomial theorem implies that this is $$n'^{(k+p^a(p-1))p^t}+n'^{(k-p^a(p-1))p^t}+kp^{a+t+1}c(n'^{(k+p^a(p-1)-1)p^t}+n'^{(k-p^a(p-1))p^t-1}) \mod p^{3a+t+v+2}$$ By assumption, we have that $p^{2a+1+v}|k$, implying that $$kp^{a+t+1}c(n'^{(k+p^a(p-1)-1)p^t}+n'^{(k-p^a(p-1))p^t-1})$$ vanishes $\mod p^{3a+t+v+2}$, and so we get that $$(nx')^{(k+p^a(p-1))p^t}+(nx')^{(k-p^a(p-1))p^t} \equiv n'^{(k+p^a(p-1))p^t}+n'^{(k-p^a(p-1))p^t} \mod p^{3a+t+v+2}$$ Hence it suffices to show that $$(nx')^{(k+p^a(p-1))p^t}+(nx')^{(k-p^a(p-1))p^t} \equiv g(n^{(k+p^a(p-1))p^t}+n^{(k-p^a(p-1))p^t}) \mod p^{3a+t+v+2}$$
We analyze $$n^{(k+p^a(p-1))p^t}(x'^{(k+p^a(p-1))p^t}-g)+n^{(k-p^a(p-1))p^t}(x'^{(k-p^a(p-1))p^t}-g)$$ Euler's Theorem implies that $x'^{(k+p^a(p-1))p^t} \equiv x'^{(k+rp^{3a+v+1}(p-1)+p^a(p-1))p^t} \mod p^{3a+t+v+2}$ for any $r$, but then notice that we can choose $r$ such that $k'|rp^{3a+v+1}(p-1)+p^a(p-1)$ (as $p$ is invertible $\mod k'$). For this $r$, we have that $\frac{(k+rp^{3a+v+1}(p-1)+p^a(p-1))p^t}{k'} \equiv 1+\frac{(rp^{3a+v+1}(p-1)+p^a(p-1))p^t}{k'} \mod d$. Since $k'|(p-1)(rp^{3a+v+1}+p^a)p^t$, we can write $k'=k_1k_2$ with $k_1|p-1$ and $k_2|(rp^{3a+v+1}+p^a)p^t$. For this $k_1$, we see that $k_1|\gcd(k,p-1)$, and so $d|\frac{p-1}{k_1}$, immediately implying that $d|\frac{(rp^{3a+v+1}(p-1)+p^a(p-1))p^t}{k'}$. Therefore, $\frac{(k+rp^{3a+v+1}(p-1)+p^a(p-1))p^t}{k'} \equiv 1 \mod d$. Letting $u=\frac{(k+rp^{3a+v+1}(p-1)+p^a(p-1))p^t}{k'}$, $$x'^{(k+p^a(p-1))p^t}-g \equiv (x'^{k'}-g)(\sum_{i=0}^{u-1} x'^{k'(u-1-i)}g^i) \mod p^{3a+t+v+2}$$ Similarly, letting $u'=\frac{(k+r'p^{3a+v+1}(p-1)-p^a(p-1))p^t}{k'}$ for some $r'$ with $k'|(r+r')(p-1)$, we have $$x'^{(k-p^a(p-1))p^t}-g \equiv (x'^{k'}-g)(\sum_{i=0}^{u'-1} x'^{k'(u'-1-i)}g^i) \mod p^{3a+t+v+2}$$ and so it remains to show that $$n^{(k+p^a(p-1))p^t}(\sum_{i=0}^{u-1} x'^{k'(u-1-i)}g^i)+n^{(k-p^a(p-1))p^t}(\sum_{i=0}^{u'-1} x'^{k'(u'-1-i)}g^i) \equiv 0 \mod p^{2a+t+v+1}$$
Now since $$\sum_{i=0}^{u-1} x'^{k'(u-1-i)}g^i=\frac{x'^{uk'}-g^u}{x'^{k'}-g}$$ lifting implies that the sum vanishes $\mod p^{v_p(u)}=p^{a+t}$, and so $$n^{(k+p^a(p-1))p^t}(\sum_{i=0}^{u-1} x'^{k'(u-1-i)}g^i) \equiv n^{kp^t}(\sum_{i=0}^{u-1} x'^{k'(u-1-i)}g^i) \mod p^{2a+2t+1}$$ by Euler's Theorem. $t \geq v$ gives this relation $\mod p^{2a+t+1+v}$, and now repeating this process with $$n^{(k-p^a(p-1))p^t}(\sum_{i=0}^{u'-1} x'^{k'(u'-1-i)}g^i)$$ gives that it suffices to show $$\sum_{i=0}^{u-1} x'^{k'(u-1-i)}g^i+\sum_{i=0}^{u'-1} x'^{k'(u'-1-i)}g^i \equiv 0 \mod p^{2a+t+v+1}$$ Write $u=p^{a+t+v}q+y$ and $u'=p^{a+t+v}q'+y'$ with $0 \leq y,y'<p^{a+t+v}$. \newline \newline 
First we claim that our expression is just $$(q+q')\sum_{i=0}^{p^{a+t+v}-1} x'^{k'(p^{a+t+v}-1-i)}g^i+g^q\sum_{i=0}^{y-1} x'^{k'(y-1-i)}g^i+g^{q'}\sum_{i=0}^{y'-1} x'^{k'(y'-1-i)}g^i \mod p^{2a+1+t+v}$$ We have that for any positive integer $c$, $$\sum_{i=cp^{a+t+v}}^{(c+1)p^{a+t+v}-1} x'^{k'(u-1-i)}g^i=\sum_{i=0}^{p^{a+t+v}-1} x'^{k'(u-1-cp^{a+t+v}-i)}g^{i+cp^{a+t+v}}$$ By lifting, we have that $x'^{k'p^{a+t+v}} \equiv g^{p^{a+t+v}} \mod p^{2a+t+v+1}$, and so $$\sum_{i=0}^{p^{a+t+v}-1} x'^{k'(u'-1-cp^{a+t+v}-i)}g^{i+cp^{a+t+v}} \equiv \sum_{i=0}^{p^{a+t+v}-1} x'^{k'(u-1-i)}g^i \mod p^{2a+t+v+1}$$ Then note that $$\sum_{i=0}^{p^{a+t+v}-1} x'^{k'(u-1-i)}g^i-x^{k'(p^{a+t+v}-1-i)}g^i=(x'^{k'(u-p^{a+t+v})}-1)\sum_{i=0}^{p^{a+t+v}-1} x^{k'(p^{a+t+v}-1-i)}g^i$$ We have that $v_p(x'^{k'(u-p^{a+t+v})}-1) \geq a+1+a+t$, and $$v_p(\sum_{i=0}^{p^{a+t+v}-1} x^{k'(p^{a+t+v}-1-i)}g^i) \geq a+t+v$$ by lifting, so $$\sum_{i=0}^{p^{a+t+v}-1} x'^{k'(u-1-i)}g^i-x^{k'(p^{a+t+v}-1-i)}g^i \equiv 0 \mod p^{3a+2t+v+1}$$ showing that $$\sum_{i=0}^{p^{a+t+v}-1} x'^{k'(u-1-i)}g^i \equiv \sum_{i=0}^{p^{a+t+v}-1} x^{k'(p^{a+t+v}-1-i)}g^i \mod p^{2a+t+v+1}$$ In all, we get that $$\sum_{i=0}^{qp^{a+t+v}-1} x'^{k'(u-1-i)}g^i \equiv q\sum_{i=0}^{p^{a+t+v}-1} x'^{k'(p^{a+t+v}-1-i)}g^i \mod p^{2a+t+v+1}$$ Now we know that $$\sum_{i=qp^{a+t+v}}^{u-1} x'^{k'(u-1-i)}g^i \equiv g^q\sum_{i=0}^{y-1} x'^{k'(y-1-i)}g^i \mod p^{2a+t+v+1}$$ by shifting since $$\sum_{i=qp^{a+t+v}}^{u-1} x'^{k'(u-1-i)}g^i \equiv \sum_{i=0}^{y-1} x'^{k'(y-1-i)}g^{i+q} \mod p^{2a+t+v+1}$$ Repeating the same arguments for $u'$ and adding finishes the proof of the claim. \newline \newline
Now we claim that this is just $$\frac{u+u'}{p^{a+t+v}}\sum_{i=0}^{p^{a+t+v}-1} x'^{k'(p^{a+t+v}-1-i)}g^i \mod p^{2a+1+t+v}$$ If $v=0$, then since $p^{a+t}|u,u'$ it immediately follows that $y=y'=0$ and $q+q'=\frac{u+u'}{p^{a+t+v}}$, proving the claim. Now suppose $v>0$, meaning that $p^{a+t+v} \nmid u,u'$. Since $0<y+y'<2p^{a+t+v}$ and $p^{a+t+v}|y+y'$, it follows that $y+y'=p^{a+t+v}$ and so $q+q'=\frac{u+u'}{p^{a+t+v}}-1$. Hence it suffices to show that $$g^q\sum_{i=0}^{y-1} x'^{k'i}g^{y-1-i}+g^{q'}\sum_{i=0}^{y'-1} x'^{k'i}g^{y'-1-i}-\sum_{i=0}^{p^{a+t+v}-1} x'^{k'i}g^{p^{a+t+v}-1-i} \equiv 0 \mod p^{2a+t+v+1}$$ Note that $$\sum_{i=0}^{y-1} x'^{k'i}g^{p^{a+t+v}-1-i}-g^qx'^{k'i}g^{y-1-i} \equiv -(g^{qp^{a+t+v}+y}-g^{p^{a+t+v}})\sum_{i=0}^{y-1} x'^{k'i}g^{-1-i} \equiv 0 \mod p^{2a+v+t+1}$$ after noting that $qp^{a+t+v}+y \equiv u \equiv 1 \mod d$.
In addition, since $$\sum_{i=y}^{p^{a+t+v}-1} x'^{k'i}g^{p^{a+t+v}-1-i}=\sum_{i=0}^{y'-1} x'^{k'(i+y)}g^{y'-1-i}$$ $$\sum_{i=y}^{p^{a+t+v}-1} x'^{k'i}g^{p^{a+t+v}-1-i}-g^{q'}\sum_{i=0}^{y'-1} x'^{k'i}g^{y'-1-i}$$ just becomes $$(x'^{k'y}-g^{q'})\sum_{i=0}^{y'-1} x'^{k'i}g^{y'-1-i} \mod p^{2a+t+v+1}$$ 
Now noting that $g^{q'} \equiv g^{q'p^{a+t+v}} \equiv x'^{k'q'p^{a+t+v}} \equiv x'^{k'(u'-y')} \mod p^{2a+t+v+1}$, it remains to show that $$(x'^{k'y}-x'^{k'(u'-y')})\sum_{i=0}^{y'-1} x'^{k'i}g^{y'-1-i} \equiv 0 \mod p^{2a+t+v+1}$$ However, $v_p(x'^{k'y}-x'^{k'(u'-y')})=v_p(x'^{k'(u'-p^{a+t+v})}-1) \geq a+1+a+t$ and $$\sum_{i=0}^{y'-1} x'^{k'i}g^{y'-1-i} \equiv 0 \mod p^{a+t}$$ so at last we conclude that $$\sum_{i=y}^{p^{a+t+v}-1} x'^{k'i}g^{p^{a+t+v}-1-i}-g^{q'}\sum_{i=0}^{y'-1} x'^{k'i}g^{y'-1-i} \equiv 0 \mod p^{2a+t+v+1}$$ completing the proof of the claim. \newline \newline Therefore, the whole problem reduces to showing that $$\frac{u+u'}{p^{a+t+v}}\sum_{i=0}^{p^{a+t+v}-1} x'^{k'(p^{a+t+v}-1-i)}g^i \equiv 0 \mod p^{2a+1+t+v}$$ Note that $u+u'$ has valuation at least $2a+1+t+v$, and so $\frac{u+u'}{p^{a+t+v}}$ has valuation at least $a+1$, so it reduces to showing that $$\sum_{i=0}^{p^{a+t+v}-1} x'^{k'(p^{a+t+v}-1-i)}g^i \equiv 0 \mod p^{a+t+v}$$ which is an immediate consequence of lifting.
\subsubsection*{Corollary 1} Let $p$ be an odd prime and $x,y \mod p$ be residues such that $x \equiv y\mu \mod p$ for some $\gcd(k,p-1)$th root of unity $\mu$. Let $\mathfrak{S}_x$ denote the multiset of residues $\mod p^{3a+t+v+2}$ attained by $n^{(k+p^a(p-1))p^t}+n^{(k-p^a(p-1))p^t}$ ($k,a,t$ as in Theorem $1$) as $n \equiv x \mod p$ varies from $1$ to $p^{a+1}$ with multiplicity and similarly for $\mathfrak{S}_y$. Then $\mathfrak{S}_x=\mathfrak{S}_y$. \newline \newline The main content of Theorem $1$ shows that if $x'^{k'} \equiv g \mod p^{a+1}$, then for any invertible $n$ and $n'$ with $n' \equiv nx' \mod p^{a+1}$, $$n'^{(k+p^a(p-1))p^t}+n'^{(k-p^a(p-1))p^t} \equiv g(n^{(k+p^a(p-1))p^t}+n^{(k-p^a(p-1))p^t}) \mod p^{3a+t+v+2}$$ Setting $g \equiv 1 \mod p^{a+1}$ and noting that $x'^{k'} \equiv 1 \mod p^{a+1}$ iff $x'^{\gcd(p-1,k)} \equiv 1 \mod p^{a+1}$ (as $\gcd(p-1,k)=\gcd(p-1,k')$, it follows that for any $\gcd(p-1,k)$th root of unity $x' \mod p^{a+1}$, $$n'^{(k+p^a(p-1))p^t}+n'^{(k-p^a(p-1))p^t} \equiv n^{(k+p^a(p-1))p^t}+n^{(k-p^a(p-1))p^t} \mod p^{3a+t+v+2}$$ Now reducing the $\gcd(p-1,k)$th roots of unity $\mod p$ shows that their reductions are $\gcd(p-1,k)$th roots of unity $\mod p$. Furthermore, if two were congruent $\mod p$, say $y,y'$, then $y^{\gcd(p-1,k)}-y'^{\gcd(p-1,k)} \equiv 0 \mod p^{a+1}$ implies that $y \equiv y' \mod p^{a+1}$ by lifting. Hence their reductions represent each $\gcd(p-1,k)$th root of unity $\mod p$ exactly once, and so let $x'$ be the unique $\gcd(p-1,k)$th root of unity $\mod p^{a+1}$ with $x' \equiv \mu \mod p$. Then note that associating $n$ and $n'$ gives a bijection between positive integers $\leq p^{a+1}$ with $n \equiv x \mod p^{a+1}$ and $n' \equiv y \mod p^{a+1}$ since the map is injective and the sets are the same size, but then since $$n'^{(k+p^a(p-1))p^t}+n'^{(k-p^a(p-1))p^t} \equiv n^{(k+p^a(p-1))p^t}+n^{(k-p^a(p-1))p^t} \mod p^{3a+t+v+2}$$ each corresponding pair yields the same element $\mod p^{3a+t+v+2}$, immediately showing $\mathfrak{S}_x=\mathfrak{S}_y$. \newline \newline We can also now prove a congruence that will be useful for the ``hard" case of the Kummer Congruences.
\subsubsection*{Corollary 2} Let $a, b, t$ be nonnegative integers with $p^{a+t}|b$ and $p-1 \nmid b$. Then we have that $$\displaystyle \sum_{n=1}^{p^{a+1}} n^b(n^{p^{a+t}(p-1)}-1)^2 \equiv 0 \mod p^{3a+t+v+2}$$ Let $S$ be the sum of all elements in $\mathfrak{S} \mod p^{3a+t+v+2}$. Then since $p-1 \nmid b$, we can find a primitive $d$th root of unity $\mod p^{3a+t+v+2}$ such that $g\mathfrak{S}=\mathfrak{S}$ for some $d>1$. If $n^{(k+p^a(p-1))p^t}+n^{(k-p^a(p-1))p^t}$ is noninvertible, then taking $\mod p$ implies $p|n$, but since $2p^{a+t} \geq 3a+2t+2$, we get that it is $0 \mod p^{3a+t+v+2}$, implying that $S \equiv gS \mod p^{3a+t+v+2}$. It suffices to show that $g \not \equiv 1 \mod p$. Suppose it is and note that $p^{3a+t+v+2}|g^d-1=(g-1)(g^{d-1}+\cdots+1)$. We have that $g^{d-1}+\cdots+1 \equiv d \mod p$, which implies the second term is coprime to $p$, and so $p^{3a+t+v+2}|g-1$, contradicting primitivity. Hence $$\sum_{n=1}^{p^{a+1}} n^{(k+p^a(p-1))p^t}+n^{(k-p^a(p-1))p^t} \equiv 0 \mod p^{3a+t+v+2}$$ and so Lemma $2$ implies that $$\sum_{n=1}^{p^{a+1}} n^{(k-p^a(p-1))p^t}(n^{p^{a+t}(p-1)}-1)^2 \equiv 0 \mod p^{3a+t+v+2}$$ Now if $p^{a+t}|b$, then $b \equiv kp^t \mod p^{a+t}(p-1)$ for some $k$ with $p^{2a+1}|k$, and so writing $b=kp^t+lp^{a+t}(p-1)$ for some integer $v$, meaning that our sum is just $$\sum_{n=1}^{p^{a+1}} n^{kp^t+lp^{a+t}(p-1)}(n^{p^{a+t}(p-1)}-1)^2-n^{kp^t-p^{a+t}(p-1)}(n^{p^{a+t}(p-1)}-1)^2 \mod p^{3a+t+v+2}$$ However, noting that $n^{|l+1|p^{a+t}(p-1)}-1$ is divisible by $n^{p^{a+t}(p-1)}-1$, the sum becomes $$\sum_{n=1}^{p^{a+1}} n^{kp^t+\min(l,-1)p^{a+t}(p-1)}f(n)(n^{p^{a+t}(p-1)}-1)^3$$ for some $f(n)$ in $\mathbb{Z}$, which vanishes $\mod p^{3a+t+v+2}$ by Euler's Theorem. \newline \newline From this, we can deduce a linearity relation that will quickly imply the Kummer congruences.
\section*{Theorem 2} Let $k$ be a positive integer not divisible by $p-1$ with $2p^{2a+1}|k$. Then we have that for any nonnegative integer $t$ and integer $r$ such that $k+p^a(p-1)r>0$, $$\displaystyle \sum_{n=1}^{p^{a+1}} n^{(k+p^a(p-1)r)p^t} \equiv r\sum_{n=1}^{p^{a+1}} n^{(k+p^a(p-1))p^t} \mod p^{3a+t+v+2}$$ \newline \newline
By Lemma $2$, we have $$\displaystyle \sum_{n=1}^{p^{a+1}} n^{kp^t} \equiv 0 \mod p^{3a+t+v+2}$$ and so the claim holds for $r=0$. \newline \newline Now  we suppose it holds for $r-1,r$. We show that it holds for $r+1$. Corollary $2$ tells us that $$\displaystyle \sum_{n=1}^{p^{a+1}} n^{{k+p^a(p-1)(r-1)}p^t}(n^{p^{a+t}(p-1)}-1)^2 \equiv 0 \mod p^{3a+t+v+2}$$ and so $$\displaystyle \sum_{n=1}^{p^{a+1}} n^{(k+p^a(p-1)(r+1))p^t} \equiv \sum_{n=1}^{p^{a+1}} 2n^{(k+p^a(p-1)r)p^t}-n^{(k+p^a(p-1)(r-1))p^t} \mod p^{3a+t+v+2}$$ By the induction hypothesis, we get $$\sum_{n=1}^{p^{a+1}} 2n^{(k+p^a(p-1)r)p^t}-n^{k+p^a(p-1)(r-1))p^t} \equiv (r+1)\sum_{n=1}^{p^{a+1}} n^{(k+p^a(p-1))p^t} \mod p^{3a+t+v+2}$$ as desired. A similar argument going backwards shows the claim for all negative $r$ with $k+p^a(p-1)r>0$ since $k+p^a(p-1)r>0$ implies $(k+p^a(p-1)r)p^t \geq 2p^{a+t} \geq 3a+2t+2$.
\section*{The proof of the congruences}
To prove the general case, we break it into three cases, the first of which is independent of the other two and the other two of which we glue together at the end. Throughout the proof, we also assume $s>r$. \newline \newline
Firstly, note that the $1-p^{s-1}$ term can always be excluded from the proofs assuming $s>r$ since then $s \geq r+p^a(p-1)>2+2^a>a+2$, meaning that $1-p^{s-1}$ always vanishes $\mod p^{a+1}$, and so since $\frac{B_s}{s}$ is a $p$-integer by Corollary $2$, we conclude that $(1-p^{s-1})\frac{B_s}{s} \equiv \frac{B_s}{s} \mod p^{a+1}$. \newline \newline
Case 1: $v_p(r)<a$ (which means $v_p(s)=v_p(r)<a$).
\newline \newline
For the sake of this case, we assume that the congruences hold $\mod p^a$ when proving them $\mod p^{a+1}$, noting that this case is empty for $a=0$. \newline \newline
Consider the sum $$\displaystyle \sum_{n=1}^{p^a} n^r(n^{p^{a-1}(p-1)}-1)^p$$ For the terms with $p$ not dividing $n$, Euler's Theorem implies that each term is divisible by $p^{ap}$ and hence $p^{2a+v_p(r)+2}$ because $ap \geq 3p \geq 2a+v_p(r)+1$. Now if $p|n$, then $p^{p^{a-1}(p-1)}|n^{p^{a-1}(p-1)}$, and so since $p^{a-1}(p-1) \geq 4(5^{a-1})>3a$, we conclude that $$\sum_{m=1}^{p^{a-1}} (pm)^r((pm)^{p^{a-1}(p-1)}-1)^p \equiv -p^r\sum_{m=1}^{p^{a-1}} m^r \mod p^{2a+v_p(r)+1}$$ We see that $$\sum_{m=1}^{p^{a-1}} m^r \equiv p^{a-1}B_r \mod p^{2a-2+v_p(r)+1}$$ by Lemma $1$, and so $$-p^r\sum_{m=1}^{p^{a-1}} m^r \equiv -p^{r+a-1}B_r \mod p^{2a+v_p(r)+1}$$ (after noting $r \geq 2$). \newline \newline Now as before, the sum expands to $$\displaystyle \sum_{n=1}^{p^a} \sum_{i=0}^{p} \binom{p}{i}(-1)^{p-i}n^{r+ip^{a-1}(p-1)}$$ Applying Lemma $1$, we get that $$\displaystyle \sum_{i=0}^{p} \binom{p}{i}(-1)^{p-i}B_{r+ip^{a-1}(p-1)} \equiv -p^{r-1}B_r \mod p^{a+v_p(r)+1}$$ and so $$\displaystyle B_{r+p^a(p-1)}-B_r+\sum_{i=1}^{p-1} \binom{p}{i}(-1)^{p-i}B_{r+ip^{a-1}(p-1)} \equiv -p^{r-1}B_r \mod p^{a+v_p(r)+1}$$ Now we have that $$\displaystyle (1-p^{r-1})\frac{B_r}{r} \equiv \frac{B_{r+ip^{a-1}(p-1)}}{r+ip^{a-1}(p-1)} \mod p^a$$ so $$B_{r+ip^{a-1}(p-1)} \equiv (1-p^{r-1})(r+ip^{a-1}(p-1))\frac{B_r}{r} \mod p^{a+v_p(r)}$$ since $v_p(r+ip^{a-1}(p-1)) \geq v_p(r)$. Now $\mod p^{a+v_p(r)+1}$, $$\displaystyle \sum_{i=1}^{p-1} \binom{p}{i}(-1)^{p-i}B_{r+ip^{a-1}(p-1)} \equiv (1-p^{r-1})\sum_{i=1}^{p-1} \binom{p}{i}(-1)^{p-i}(r+ip^{a-1}(p-1))\frac{B_r}{r}$$ and noting $$\displaystyle \sum_{i=1}^{p-1} \binom{p}{i}(-1)^{p-i}i=-p$$ gives $-(1-p^{r-1})p^a(p-1)\frac{B_r}{r}$. We have that $v_p(p^{a+r-1}\frac{B_r}{r}) \geq a+r-1 \geq a+v_p(r)+1$, invoking Lemma $3$ and noting $r \geq v_p(r)+2$, and so we conclude that $-(1-p^{r-1})p^a(p-1)\frac{B_r}{r} \equiv -p^a(p-1)\frac{B_r}{r} \mod p^{a+v_p(r)+1}$. Putting everything together gives $B_{r+p^a(p-1)} \equiv (r+p^a(p-1))\frac{B_r}{r}-p^{r-1}B_r \mod p^{a+v_p(r)+1}$. Lemma $3$ again gives $v_p(p^{a+r-1}B_r) \geq a+r-1+v_p(r)$, and so since $r \geq 2$, we see that $p^a(p-1)(p^{r-1}B_r)$ vanishes $\mod p^{a+v_p(r)+1}$, and so $$B_{r+p^a(p-1)} \equiv (r+p^a(p-1))\frac{B_r}{r}-(r+p^a(p-1))p^{r-1}\frac{B_r}{r} \mod p^{a+v_p(r)+1}$$ Finally, since $v_p(r)=v_p(r+p^a(p-1))$, we conclude that $$\frac{B_{r+p^a(p-1)}}{r+p^a(p-1)} \equiv (1-p^{r-1})\frac{B_r}{r} \mod p^{a+1}$$ finishing this case by iterating this congruence if necessary. \newline \newline
To handle the next case, we must make use of Theorem $2$ instead of just using direct induction. This is because if we try to repeat the same argument that we did for the case of $v_p(r)<a$, we will run into the problem where the binomial coefficients only give one extra power of $p$ when we actually need $2$ (one to get up to $a$ and one because the powers separating $r$ are larger than those we would have in the sum) \newline \newline Note that the $1-p^{r-1}$ may now also be excluded from now on since $r \geq 2p^a \geq 2^{a+1} \geq a+2$ and the same arguments used for $s$. \newline \newline
Case 2: $v_p(r)=v_p(s) \geq a$. \newline \newline
Theorem $2$ immediately implies that $\mod p^{4a+t+2}$, $$(k+p^a(p-1))\displaystyle \sum_{n=1}^{p^{a+1}} n^{(k+p^a(p-1)b)p^t} \equiv (k+bp^a(p-1)) \sum_{n=1}^{p^{a+1}} n^{(k+p^a(p-1))p^t}$$ Now, by Lemma $1$, we have that $$(k+p^a(p-1))B_{(k+bp^a(p-1))p^t} \equiv (k+bp^a(p-1))B_{(k+p^a(p-1))p^t} \mod p^{3a+t+1}$$ and so if $r$ is not a multiple of $p$, then dividing by $p^t(k+p^a(p-1))(k+bp^a(p-1))$ (of valuation $2a+t$) immediately gives the theorem. 
\newline \newline Case 3: Branching \newline \newline
The key now is to consider the object $$\displaystyle \sum_{n=1}^{p^{a+1}} n^{kp^t}(n^{p^{a+t-1}(p-1)}-1)^p$$ which vanishes $\mod p^{p(a+t)}$ and hence $\mod p^{3a+t+2}$. By the binomial theorem, we have that this is just \newline $$\displaystyle \sum_{n=1}^{p^{a+1}} n^{(k+p^a(p-1))p^t}-n^{kp^t}+\sum_{n=1}^{p^{a+1}} \sum_{i=1}^{p-1} \binom{p}{i}(-1)^{p-i}n^{kp^t+ip^{a+t-1}(p-1)}$$ The double sum is just $$ \sum_{i=1}^{p-1} \binom{p}{i}(-1)^{p-i}(\frac{k}{p^a}+i) \sum_{n=1}^{p^{a+1}} n^{(k+p^a(p-1))p^{t-1}} \mod p^{3a+t+2}$$ by Theorem $2$, so the whole sum just reduces to $$-p\sum_{n=1}^{p^{a+1}} n^{(k+p^a(p-1))p^{t-1}} \mod p^{3a+t+2}$$ This implies that $$\displaystyle \sum_{n=1}^{p^{a+1}}n^{(k+p^a(p-1))p^t} \equiv p\sum_{n=1}^{p^{a+1}} n^{(k+p^a(p-1))p^{t-1}} \mod p^{3a+t+2}$$ Using Lemma $1$ once again, we get that $$B_{(k+p^a(p-1))p^t} \equiv pB_{(k+p^a(p-1)r)p^{t-1}} \mod p^{2a+t+1}$$ Finally, dividing by $(k+p^a(p-1))p^t$ (of valuation $a+t$) immediately gives $$\frac{B_{(k+p^a(p-1))p^t}}{(k+p^a(p-1))p^t} \equiv  \frac{B_{(k+p^a(p-1))p^{t-1}}}{(k+p^a(p-1))p^{t-1}} \mod p^{a+1}$$ completing the branching step and hence the proof of the full Kummer congruences. \newline \newline
Thus we've completed our goal of proving the Kummer Congruences by understanding the local behavior and interpreting that in the context of a global sum. \newline \newline Notice that the set $\mathfrak{S}$ comes equipped with a $(\mathbb{Z}/p^{3a+t+v+2}\mathbb{Z})^*$ action and that the content of Theorem $3$ is that the subgroup $\mu_d$ of $d$th roots of unity is contained in $\text{Stab}(\mathfrak{S})$ under this action. We now classify $\text{Stab}(\mathfrak{S})$ more completely. We need some lemmas. \newline \newline Before starting the proof, we make a quick note about $p$-adic valuations and polynomials. For a polynomial $g(x)$, we say that for a residue class $s \mod p^k$, $v_p(g(s)) \geq l$ iff $g(s') \equiv 0 \mod p^l$ for every positive integer $s'$ in the residue class $s \mod p^k$ (in general, we can say that $v_p(g(s))=\min_{s' \equiv s \mod p^k} v_p(g(s'))$.
\subsection*{Lemma 4} Let $p$ be an odd prime and set $f(x)=x^{(k+p^a(p-1))p^t}+x^{(k-p^a(p-1))p^t}$ as polynomials. Then for any positive integers $m,n$ with $n$ coprime to $p$, we have that $v_p(f^{(m)}(n)) \geq 2a+t+v$. Furthermore, if $m=1$, we have that $v_p(f^{(m)}(n)) \geq 2a+t+1+v$, with equality always holding if $v<t$. For $m=2$, we have that $v_p(f^{(m)}(n))=2a+2t$ if $v=t$ and that $v_p(f^{(m)}(n) \geq 2a+t+v+1$ if $v<t$. Finally, for $m=3$ and $v=t$, we have that $v_p(f^{(m)}(n))=2a+2t$ unless $p=3$, in which case it is $\geq 2a+2t+1$, and for $v<t$, $v_p(f^{(m)}(n)) \geq 2a+t+v+1$. \newline \newline Expanding $f^{(m)}(n)$, we get $$n^{(k+p^a(p-1))p^t-m}\prod_{i=0}^{m-1} ((k+p^a(p-1))p^t-i)+n^{(k-p^a(p-1))p^t-m}\prod_{i=0}^{m-1} ((k-p^a(p-1))p^t-i)$$ Since $m$ is a positive integer, each product has valuation at least $a+t$, and so Euler's Theorem implies that this is just $$ n^{kp^t-m}(\prod_{i=0}^{m-1} ((k+p^a(p-1))p^t-i)+\prod_{i=0}^{m-1} ((k-p^a(p-1))p^t-i)) \mod p^{2a+2t+1}$$ Then note that the sum of products is just $$ p^{a+t}(p-1)(\prod_{i=1}^{m-1} (p^{a+t}(p-1)-i)+(-1)^m\prod_{i=1}^{m-1} (p^{a+t}(p-1)+i)) \mod p^{2a+t+1+v}$$ Plugging in $m=1$ immediately gives the inequality in the second statement. For the equality case, plugging in $m=1$ into the expression $$n^{kp^t-m}\prod_{i=0}^{m-1} ((k+p^a(p-1))p^t-i)+\prod_{i=0}^{m-1} ((k-p^a(p-1))p^t-i))$$ gives $f'(n) \equiv 2n^{kp^t-m}kp^t \mod p^{2a+2t+1}$, and if $v_p(k)<2a+1+t$, then $v_p(kp^t)<2a+2t+1$, so the valuation of $f'(n)$ is exactly $v_p(k)+t=2a+t+1+v$ in this case. \newline \newline For the case of $m=2$, we have that the expression is $2p^{2a+2t}(p-1) \mod p^{2a+t+v+1}$. If $t=v$, then we get that $2p^{2a+2t}(p-1) \mod p^{2a+2t+1}$, and so $v_p(f''(n))=2a+2t$. If $v<t$, note that $2a+2t \geq 2a+t+v+1$, and so $2p^{2a+2t}(p-1) \equiv 0 \mod p^{2a+t+v+1}$, as desired. \newline \newline Plugging in $m=3$, $p^{a+t}(p-1)((p^{a+t}(p-1)-1)(p^{a+t}(p-1)-2)-(p^{a+t}(p-1)+1)(p^{a+t}(p-1)+2)) \equiv -6p^{2(a+t)}(p-1)^2 \mod p^{2a+t+v+1}$. If $v=t$, then we conclude that $f'''(n)$ has valuation $2a+2t$ for $p \neq 3$ and vanishes $\mod p^{2a+2t+1}$ for $p=3$. If $v<t$, then $f'''(n)$ automatically vanishes $\mod p^{2a+t+v+1}$ once again. \newline \newline For the first, it suffices to show that $$\prod_{i=1}^{m-1} (p^{a+t}(p-1)-i)+(-1)^m\prod_{i=1}^{m-1} (p^{a+t}(p-1)+i) \equiv 0 \mod p^{a+v}$$ Reducing $\mod p^{a+v}$ gives $$(-1)^{m-1}(m-1)!+(-1)^m(m-1)! \equiv 0 \mod p^{a+v}$$ as desired.
\subsubsection*{Corollary 3} Let $r$ be a nonnegative integer and $k \geq r$ be a positive integer. Then for any integer $s$ coprime to $p$, $f(s+p^kx) \equiv f(s)+f'(s)p^kx \mod p^{2a+t+v+2k}$. Furthermore, if $v<t$, then we have the stronger $f(s+p^kx) \equiv f(s)+f'(s)p^kx \mod p^{2a+t+v+2k+1}$\newline \newline Let $\text{deg}(f)=w$. We have that for any positive integer $k$, $$f(s+p^kx)=\sum_{n=0}^{w} f^{(n)}(s)\frac{(p^kx)^n}{n!}$$ Lemma $4$ implies that $f^{(n)}(s) \geq 2a+t+v$, and so for $n>1$, we have that $v_p(f^{(n)}(s)\frac{(p^kx)^n}{n!}) \geq 2a+t+v+2k$, and so $$f(s+p^kx) \equiv f(s)+f'(s)p^kx \mod p^{2a+t+v+2k}$$ Now suppose that $v<t$. Then we claim that $v_p(f^{(n)}(s)\frac{(p^kx)^n}{n!}) \geq 2a+t+v+2k+1$ for all $n \geq 2$. For $n=2$ or $3$, Lemma $4$ implies that $v_p(f^{(n)}(s)) \geq 2a+t+v+1$, and so $v_p(f^{(n)}(s)\frac{(p^kx)^n}{n!}) \geq 2a+t+v+1+2k$ for $n=2$ and $2a+t+v+1+3k-1 \geq 2a+t+v+1+2k$ for $n=3$, proving the claim in these cases. For $n>3$, we have that $v_p(f^{(n)}(s)\frac{(p^kx)^n}{n!}) \geq 2a+t+v+4k-1 \geq 2a+t+v+2k+1$, finishing the proof in all cases and concluding the case $v<t$.
\subsection*{Lemma 5} Suppose $v=t$. Then for each nonnegative integer $s$ with $0 \leq s \leq a$, there exist exactly $p^{a-s}(p-1)$ positive integers $u \leq p^{a+1}$ with $u$ coprime to $p$ and $v_p(f'(u)) \geq 2a+2t+s+1$. \newline \newline We show by induction that there are $p-1$ solutions $u \mod p^{s+1}$ to $f'(x) \equiv 0 \mod p^{2a+2t+s+1}$ for any nonnegative integer $s$. The base case of $s=0$ is immediate from Lemma $4$. Suppose that the statement holds for $s-1$ and let $u'$ be a solution to $f'(x) \equiv 0 \mod p^{2a+2t+s}$. We show that $u'$ can be lifted to $u \mod p^{s+1}$ with $f'(u) \equiv 0 \mod p^{2a+2t+s+1}$. By Taylor series expansion, we have that $$f'(u'+p^sx)=\sum_{n=0}^{w-1} f^{(n+1)}(u')\frac{(p^sx)^n}{n!}$$ Taking out the first two terms gives $$f'(u')+f''(u')p^sx+\sum_{n=2}^{w-1} f^{(n+1)}(u')\frac{(p^sx)^n}{n!}$$ Lemma $4$ implies that for any $n$, $f^{(n+1)}(u')$ vanishes $\mod p^{2a+2t}$, and since for $n>1$ we have that $v_p(\frac{(p^sx)^n}{n!}) \geq 2s \geq s+1$ since $s \geq 1$ by assumption, $f^{(n+1)}(u')\frac{(p^sx)^n}{n!}$ vanishes $\mod p^{2a+2t+s+1}$. Hence $f'(u'+p^sx) \equiv f'(u')+f''(u')p^sx \mod p^{2a+2t+s+1}$, immediately implying that there is a unique lift of $u'$. Furthermore, for any $u' \mod p^s$ with $f'(u') \not \equiv 0 \mod p^{2a+2t+s}$, we have that $$f'(u'+p^sx) \equiv f'(u')+f''(u')p^sx \mod p^{2a+2t+s+1}$$ but then $v_p(f'(u'))<v_p(f''(u')p^s)$, and so there is no lift giving a solution here, meaning that all solutions $\mod p^{s+1}$ are lifts of solutions $\mod p^s$. With this proven, then note that there are $p^{a-s}(p-1)$ positive integers $k\leq p^{a+1}$ with $k \equiv u \mod p^{a+1}$, completing the proof.
\section*{Theorem 3} Given the conditions of Theorem $1$, $\text{Stab}(\mathfrak{S})=\mu_n$ where $n=dp^a$ if $v<t$ and $n=d$ if $v=t$. \newline \newline Since $\text{Stab}(\mathfrak{S})$ is defined by the group action of $(\mathbb{Z}/p^{3a+t+v+2}\mathbb{Z})^*$ on $\mathfrak{S}$, we have that $\text{Stab}(\mathfrak{S})$ is a subgroup of $(\mathbb{Z}/p^{3a+t+v+2}\mathbb{Z})^*$, and so is in particular cyclic. Letting $P$ be the $p$-Sylow subgroup of $\text{Stab}(\mathfrak{S})$, we can write $\text{Stab}(\mathfrak{S})=PG$ for some  abelian group $G$ of order coprime to $p$ (where the product is isomorphic to a direct one). $G$, then, consists of all elements of order coprime to $p$ and hence is of order dividing $p-1$. \newline \newline Since $p$ is odd, the subgroup of such elements is automatically $\mu_{p-1}$ since $(\mathbb{Z}/p^{3a+t+v+2}\mathbb{Z})^*$ is cyclic. Now note that if $g \in \text{Stab}(\mathfrak{S})$, its reduction $g'$ $\mod p$ is automatically in $\text{Stab}(\mathfrak{S'})$, where $S'$ is $S$ defined $\mod p$. We have that $$n^{(k+p^a(p-1))p^t}+n^{(k-p^a(p-1))p^t} \equiv 2n^{k'} \mod p$$ and so $\mathfrak{S'}$ consists of elements of the form $2n^{\gcd(k,p-1)}$, which is a coset of $\mu_d'$. Hence $\text{Stab}(\mathfrak{S'})=\mu_d'$, meaning $g' \in \mu_d'$. However since $\text{ord}_p(g)|\text{ord}_{p^{3a+t+v+2}}(g)$ and $\text{ord}_{p^{3a+t+v+2}}(g)$ is coprime to $p$, lifting implies that $\text{ord}_p(g)=\text{ord}_{p^{3a+t+v+2}}(g)$ and hence $g \in \mu_d$. Theorem $2$ then implies that $\mu_d \subset \text{Stab}(\mathfrak{S})$ and in particular that $\mu_d \subset G$, so we conclude that $G=\mu_d$. \newline \newline It remains to analyze $P$. First we find a characterization for when $P$ has order $\leq p^r$ for some nonnegative integer $r$. Recalling that it is cyclic, the condition is equivalent to there being no element of order $p^{r+1}$. If $g$ were some element of order $p^{r+1}$, then it would necessarily be $1 \mod p^{3a+t+v+1-r}$, so for any element $s \in \mathfrak{S}$, the orbit of $s$ under the subgroup generated by $g$ would be the set consisting of $s+p^{3a+1+t+v-r}x \mod p^{3a+t+v+2}$ for $0 \leq x \leq p^{r+1}-1$. Hence each element $s'$ with $s' \equiv s \mod p^{3a+1+t+v-r}$ must appear the same number of times in $\mathfrak{S}$. \newline \newline Now we introduce some terminology. Given nonnegative integers $i \geq k$ and a particular residue $s \mod p^k$, we define the values of $g(s) \mod p^{2a+t+1+v+i}$ on the $p^{i-k}$ lifts of $s \mod p^i$ to be the branches of $g(s)$ of level $i$. Furthermore, we define the layer $i$ of such an $s$ to be the multiset of all branches of level $i$ (taken $\mod p^{2a+t+1+v+i}$). For $t>j$, we call a multiset $\mathfrak{T} \mod p^t$ $j$-balanced if for any element $s$ in $\mathfrak{T}$, the $p^j$ lifts of its reduction $\mod p^{t-j}$ appear an equal number of times in $\mathfrak{T}$. Furthermore, we say the $j$-character of a multiset is good if it is $j$-balanced and bad otherwise. Note that deleting a good multisubset never changes the $j$-character of the multiset. \newline \newline Under these definitions, we have that $P$ has size $\leq p^r$ iff $\mathfrak{S}$ is not $r+1$-balanced. Now note that $\mathfrak{S}$ is the disjoint union of $\mathfrak{S}_x$ as $x$ ranges over the invertible residues $\mod p$. Let $s_1,...,s_d$ be a set of coset representatives of $\mu_{\gcd(p-1,k)}$ in $(\mathbb{Z}/p\mathbb{Z})^*$ By Corollary $3$, $\mathfrak{S}$ is just $\gcd(p-1,k)$ copies of the multiset $\coprod_{i=1}^{d} \mathfrak{S}_{s_i}$. Let $\mathfrak{V}=\coprod_{i=1}^{d} \mathfrak{S}_{s_i}$. It immediately follows that $\mathfrak{V}$ is $r+1$-balanced iff $\mathfrak{S}$ is. Noting again that $$n^{(k+p^a(p-1))p^t}+n^{(k-p^a(p-1))p^t} \equiv 2n^{k'} \mod p$$ if $x \not \equiv y\mu \mod p$ for some $\gcd(p-1,k)$th root of unity $\mu$, it follows that $\mathfrak{S}_x$ and $\mathfrak{S}_y$ are disjoint even upon restricting $\mod p$. Hence for any $s \in \mathfrak{S}_x$, no element $s'$ with $s' \equiv s \mod p^{3a+1+t-r}$ can appear in $\mathfrak{S}_y$, and so $\mathfrak{V}$ is $r+1$-balanced iff for each $x$ and $s \in \mathfrak{S}_x$, all $s'$ appear the same number of times in $\mathfrak{S}_x$, which is precisely equivalent to $\mathfrak{S}_x$ being $r+1$-balanced for each $x$. Finally, note that if $\mathfrak{S}_x$ is $r+1$-balanced, then it follows that $\mathfrak{S}'_x$ is $r+1-q$-balanced, where $\mathfrak{S}'_x$ is the reduction of $\mathfrak{S}_x \mod p^{3a+t+v+2-q}$. \newline \newline 
Case 1: $v<t$ \newline \newline In this case, we show that $P$ is cyclic of order $p^a$. To show that $P$ has size $p^a$, we show that each $\mathfrak{S}_x$ is $a$-balanced but not $a+1$-balanced. To show that $\mathfrak{S}_x$ is not $a+1$-balanced, note that taking $k=1$ in Corollary $3$ gives $f(s+px) \equiv f(s)+f'(s)px \equiv f(s) \mod p^{2a+t+v+2}$, so the reduction of $\mathfrak{S}_x \mod p^{2a+t+v+2}$ fails to be $1$-balanced, implying that $\mathfrak{S}_x$ fails to be $a+1$-balanced for each $x$. \newline \newline Now we show that each $\mathfrak{S}_x$ is $a$-balanced. Consider the sequence of positive integers given by $a_n=\lceil \frac{a+1}{2^n} \rceil$ for all $n \geq 0$. Noting that $2\lceil x \rceil \geq \lceil 2x \rceil$, we see that $2a_n \geq a_{n-1}$, and so we have that $f(s+p^{a_n}x) \equiv f(s)+f'(s)p^{a_n}x \mod p^{2a+t+v+a_{n-1}+1}$ for each positive integer $n$ by the stronger version of Corollary $3$. We now show that $\mathfrak{S}_x$ is $a_0-a_k=a+1-a_k$ balanced for each positive integer $k$. We pick a residue $s \mod p^{a_k}$ and count how many times a fixed lift $u$ of $f(s) \mod p^{2a+t+v+a_k+1}$ appears in the $a+1$ layer of $s$. To do this, it suffices to find the number of $s' \mod p^{a+1}$ with $s' \equiv s \mod p^{a_k}$ and $f(s') \equiv u \mod p^{3a+t+2}$. We show by induction that there is a unique lift $s'$ of $s \mod p^{a_r}$ that $u \equiv f(s') \mod p^{2a+t+v+a_r+1}$ for all $0 \leq r \leq n$. The base case of $r=k$ is immediate. Now suppose that it holds for $n$ and we show it for $n-1$. Let $s_n$ be the unique lift of $s \mod p^{a_n}$ such that $u \equiv f(s_n) \mod p^{2a+t+v+a_n+1}$. Then we have that $f(s_n+p^{a_n}x) \equiv f(s_n)+f'(s_n)p^{a_n}x \mod p^{2a+t+v+a_{n-1}+1}$. As $v_p(f'(s_n)p^{a_n})=2a+t+v+a_n+1$, this uniquely determines $x \mod p^{a_{n-1}-a_n}$ and hence uniquely determines $s_{n-1}$. Therefore, $s_0$ is uniquely determined $\mod p^{a+1}$, showing that there are exactly $1$ lift of $s \mod p^{a+1}$ corresponding to a given $u$. Hence $\mathfrak{S}_x$ is $a_0-a_k=a+1-a_k$ balanced. Taking $k$ large enough implies that $\mathfrak{S}_x$ is $a$-balanced. \newline \newline
Case 2: $t=v$ \newline \newline From Corollary $3$, given that $v_p(f'(s))=2a+2t+1+u$ for some nonnegative integer $u$, $f(s+p^kx) \equiv f(s) \mod p^{2a+2t+k+u+1}$ for all $k$ with $k \geq u+1$, showing that for $s'$ the reduction of $s \mod p^k$, $f(s')$ is constant on its $k+u$ branches for such $k$. However, we also have $f(s+p^kx) \equiv f(s)+f'(s)p^kx \mod p^{2a+2t+k+u+2}$ for $k \geq u+2$, so for such $k$, while the $k+u+1$ branches of $f(s')$ of a given $k+u$ branch of $f(s')$ are constant, the $k+u+1$ branches of different $k+u$ branches are distinct, hence comprising all lifts $\mod p^{2a+k+u+1}$ of the $k+u$ branch of $f(s')$ once each, implying that the $k+u+1$ layer of $s'$ is $1$-balanced. \newline \newline Now suppose that for some positive integer $u \leq \lfloor \frac{a}{2} \rfloor$, we have an $s \mod p^{u+1}$ with $v_p(f'(s))<2a+2t+u+1$. Then the first condition and the proof of Lemma $5$ imply that for any lift $s' \mod p^k$ of $s$, we have $v_p(f'(s'))=v_p(f'(s))$. Letting $v_p(f'(s))=2a+2t+1+w$, $w<u$ implies that $w \leq \lfloor \frac{a}{2} \rfloor-1$. In particular, selecting $k=3a+2t+1-v_p(f'(s'))$, the inequality on $w$ implies that $k=a-w \geq w+2$, so the above applies and hence the $a+1$ layer lying above $s' \mod p^k$ is $1$-balanced, noting that $v_p(f'(s))=v_p(f'(s'))$, and so the $a+1$ layer of $s$, which is the disjoint union of the $a+1$ layers of its lifts $\mod p^k$, is also balanced. Looking at each individual $\mathfrak{S}_x$, Lemma $5$ now allows us to remove the $a+1$ layers above all but $1$ element $\mod p^{\lfloor \frac{a}{2} \rfloor+1}$ without changing the $1$-quality of $\mathfrak{S}$, so it remains to show that the $a+1$ layer above any (though we show each) exceptional element $\mod p^{\lfloor \frac{a}{2} \rfloor+1}$ isn't $1$-balanced. By Corollary $3$, we have that for any such exceptional $s$, $f(s+p^{\lfloor \frac{a}{2} \rfloor+1}x) \equiv f(s)+f'(s)p^{\lfloor \frac{a}{2} \rfloor+1}x \mod p^{2a+2t+2\lfloor \frac{a}{2} \rfloor+2}$, but then since $f'(s) \equiv 0 \mod p^{2a+2t+\lfloor \frac{a}{2} \rfloor+1}$, we just get $f(s+p^{\lfloor \frac{a}{2} \rfloor+1}x) \equiv f(s) \mod p^{2a+2t+2\lfloor \frac{a}{2} \rfloor+2}$. \newline \newline If $a$ is even, $2a+2t+2\lfloor \frac{a}{2} \rfloor+2=3a+2t+2$, and so all $a+1$ branches of each exceptional $s$ are the same, immediately implying that no $\mathfrak{S}_x$ is $1$-balanced. \newline \newline It remains to handle the case of $a$ odd. Taking an exceptional element $s \mod p^{\frac{a+1}{2}}$, it once again suffices to show that its $a+1$ layer is not $1$-balanced. By the proof of Lemma $5$, we can choose a positive integer $s'$ with reduction $s \mod p^{\frac{a+1}{2}}$ such that $f'(s') \equiv 0 \mod p^{\frac{5a+3}{2}+2t}$. By Taylor series expansion, we have that $$f(s'+p^\frac{a+1}{2}x)=\sum_{n=0}^{w} f^{(n)}(s')\frac{(p^{\frac{a+1}{2}}x)^n}{n!}$$ For $n>2$, we claim that $v_p(f^{(n)}(s')\frac{(p^{\frac{a+1}{2}}x)^n}{n!}) \geq 3a+2t+2$. For $n>3$, note that since Lemma $4$ implies $v_p(f^{(n)}(s')) \geq 2a+2t$, $v_p(f^{(n)}(s')\frac{(p^{\frac{a+1}{2}}x)^n}{n!}) \geq 4a+2t+1 \geq 3a+2t+2$ since $a \geq 1$. For $n=3$ and $p>3$, the inequality $v_p(f^{(n)}(s') \geq 2a+2t$ again gives $v_p(f^{(n)}(s')\frac{(p^{\frac{a+1}{2}}x)^n}{n!}) \geq 2a+2t+3(\frac{a+1}{2}) \geq 3a+2t+2$. Finally, for $n=3, p=3$, Lemma $4$ implies $v_p(f^{(n)}(s') \geq 2a+2t+1$, and so we still recover $v_p(f^{(n)}(s')\frac{(p^{\frac{a+1}{2}}x)^n}{n!}) \geq 2a+2t+3(\frac{a+1}{2}) \geq 3a+2t+2$, proving the claim. \newline \newline Hence we conclude that $$f(s'+p^\frac{a+1}{2}x)= f(s')+f'(s')p^\frac{a+1}{2}x+f''(s')\frac{p^{a+1}x^2}{2} \mod p^{3a+2t+2}$$ By construction of $s'$, $v_p(f'(s')p^\frac{a+1}{2}x) \geq \frac{5a+3}{2}+2t+\frac{a+1}{2}=3a+2t+2$, and so $$f(s'+p^\frac{a+1}{2}x)= f(s')+f''(s')\frac{p^{a+1}x^2}{2} \mod p^{3a+2t+2}$$ It suffices to show that $f''(s')\frac{p^{a+1}x^2}{2}$ cannot attain all residues of the form $p^{3a+2t+1}k \mod p^{3a+2t+2}$, which is an immediate consequence of the fact that the map $x^2 \mod p$ fails to be surjective, completing the proof.

\section*{References}
[1] Kubota, Tomio; Leopoldt, Heinrich-Wolfgang (1964), ``Eine p-adische Theorie der Zetawerte. I. Einführung der p-adischen Dirichletschen L-Funktionen," Journal für die reine und angewandte Mathematik, 214/215: 328–339 \newline \newline
[2] Kummer, Ernst Eduard (1851), ``Über eine allgemeine Eigenschaft der rationalen Entwicklungscoëfficienten einer bestimmten Gattung analytischer Functionen," Journal für die Reine und Angewandte Mathematik, 41: 368–372 \newline \newline
[3] G.F. Voronoi, On Bernoulli numbers, Comm. Charkov Math. Soc. 2 (1890) 129-148 (Russian) \newline \newline
\fontfamily{lmss}\selectfont {Department of Mathematics, California Institute of Technology, Pasadena, CA, 91125} \newline \fontfamily{lmss}\selectfont E-mail address: sgoodman@caltech.edu
\end{document}